\documentclass[a4paper]{article}

\usepackage{amsmath,amsfonts,latexsym,theorem}
\nonstopmode

\numberwithin{equation}{section}

\newtheorem{Theorem}{Theorem}[section]
\newtheorem{Definition}{Definition}[section]
\newtheorem{Proposition}{Proposition}[section]
\newtheorem{Lemma}{Lemma}[section]

\newenvironment{Proofc}[1]{\smallskip\par\noindent\textsc{#1}\quad}
  {\hfill$\Box$\bigskip\par}
\newenvironment{Proof}{\begin{Proofc}{Proof}}{\end{Proofc}}

\newtheorem{Remark}{Remark}[section]

\def\a{\alpha}

\def\d{\delta}

\def\G{\Gamma}

\def\pd{\partial}

\newcommand{\p}{\partial}

\newcommand{\R}{{\mathbb R}}
\newcommand{\N}{{\mathbb N}}

\def\pd{\partial}

\begin{document}
\title{A comparison among   various notions of viscosity solutions  for Hamilton-Jacobi equations
on networks}

\author{Fabio Camilli\footnotemark[1] \and Claudio Marchi \footnotemark[2]}

\date{version: \today}
\maketitle

\footnotetext[1]{Dip. di Scienze di Base e Applicate per l'Ingegneria,  ``Sapienza" Universit{\`a}  di Roma, via Scarpa 16,
 00161 Roma, Italy, ({\tt e-mail:camilli@dmmm.uniroma1.it})}
\footnotetext[2]{Dip. di Matematica, Universit\`a di Padova, via Trieste 63, 35121 Padova, Italy ({\tt marchi@math.unipd.it}).}

\begin{abstract}

Three definitions of viscosity solutions for Hamilton-Jacobi equations on networks recently appeared in literature (\cite{acct,imz,sc}). Being motivated by various  applications, they appear to be  considerably different. Aim of this note is to establish their equivalence.
%
\end{abstract}
\begin{description}
\item [{\bf MSC 2000}:] 35R02, 35F21, 35D40.
\item [\textbf{Keywords}:] Hamilton-Jacobi equation;  viscosity solution; network.
\end{description}
\section{Introduction}
The theory of viscosity solutions (see~\cite{bcd} for an overview) has been extensively studied and refined by many authors, and, among the numerous contributions in the literature, one can find several adaptations to very different settings. In the recent time, there is an increasing interest in the study of of nonlinear PDEs on networks and, concerning Hamilton-Jacobi equations,   three different notions of viscosity solution   have been introduced (\cite{acct}, \cite{imz}, \cite{sc}).\par
A major task in the theory of PDEs on network   is to establish the correct transition conditions the  solutions
are subjected to at vertices. It is easy to see that classical transition conditions such the Kirchhoff condition,
based on the divergence structure of the problem,   are not adequate to characterize the expected viscosity solution of the equation. Hence  in all the three approaches  the equation is   considered  also at the vertices. On the other side, since the three papers are motivated by different applications (respectively, a control problem constrained to a network in \cite{acct}, the study of traffic flow at a junction in \cite{imz} and Eikonal equations and  distance functions on  networks in \cite{sc}), they differ  for the assumptions made on the Hamiltonians and mainly for the definitions of viscosity solution at the vertices (while inside the edges all the definitions coincide with the classical one). \par
Nevertheless, since    all the definitions  give  existence and uniqueness of the solution, it is worth to compare them. In this paper, we show that, at   least when restated in a common framework, the three definitions are equivalent.
Obviously, imposing common assumptions to the problems would restrict their generality; for example, the comparison of the definitions in \cite{imz} and \cite{sc} should require that the Hamiltonian depends only on the gradient (i.e., it is independent of the state variable and of the edge) and it is strictly convex.
However, for the sake of generality, in this paper we shall keep assumptions as weak as possible, often even weaker of the original ones.

This paper is organized as follows: in section \ref{sect:setting} we introduce the network and the Hamilton-Jacobi equation we consider on it. Section \ref{sect:definition} contains the three definition of solution. Section \ref{sect:acctimz} (respectively, section \ref{sect:imzcs}) is devoted to establish the equivalence between the definitions in \cite{acct} (resp., in \cite{sc}) and the one in \cite{imz}.
%
%
%
%
\section{Setting of the problem}\label{sect:setting}
We consider a planar star-shaped network $\G$
composed  of a transition vertex $v$ and  of a finite number of straight edges $e_j$, $j\in J\equiv\{1,\dots,N\}$, i.e.
\begin{equation}\label{star}
\G= \{v\}\cup \bigcup_{j\in J} e_j \subset \R^2,\quad v=(0,0), \quad e_j=(0,1)\eta_j
\end{equation}
where $(\eta_j)_{j\in J}$ is a set of unit vectors in $\R^2$ with $\eta_j\not=\eta_k$ if $j\not =k$.
For each edge $e_j$, we fix a  parametrization $\pi_j:[0,l_j]\to\R^2$ such that  $e_j=\pi_j((0,l_j))$ and $\bar e_j=\pi_j([0,l_j])$.
Moreover we assume: $v=\pi_j(0)$ for any $j\in J$; in this way, we fix an orientation of $e_j$.

\vskip3mm

Consider a function $u:\bar \G\to \R$; for $j\in J$, we denote by   $u^j:[0,l_j]\to \R$ the restriction of $u$ to $\bar e_j$, i.e.  $u^j(y)=u(x)$ for $y=\pi_j^{-1}(x)\in [0,l_j]$.
We say that $u$ is continuous (resp., upper or lower semi-continuous) when it is so with respect to the topology induced on $\bar \G$ from $\R^2$ and we shall write $u\in C(\bar \G)$ (resp., $u\in USC(\bar \G)$ or $u\in LSC(\bar \G)$).
As in \cite{imz, sc}, we consider derivate with respect to the parametrization
\begin{equation}\label{deriv}
D_ju(x):= \frac{d  u^j}{d y }  (y)\qquad \text{for $y=\pi^{-1}_j(x)$}.
\end{equation}
We denote by $\pd^+_ju(v)$ the super-differential of $u$ at $v$ along the edge $e_j$, i.e.
\[\pd^+_ju(v):=\{p\in\R:\, u^j(y)\le u^j(0)+py+o(y)\quad\text{for  $x\in e_j\to v$, $y=\pi^{-1}_j(x)$}\}\]
and we set $\pd^-_ju(v)=-\pd^+_j(-u(v))$.

\vskip3mm

We consider the Hamilton-Jacobi equation
\begin{equation}\label{HJ}
 u+ H(x, Du )=0,\qquad x\in\G
\end{equation}
where $u:\G\to \R$ and $H:\G\times \R\to \R$.
In the following  $H^j:[0,l_j]\times \R\to\R$  denotes  the restriction of $H$ to $\bar e_j$; throughout this paper, we shall assume
\begin{equation}\label{ipofisse}
H^j\in C^0([0,l_j]\times \R), \qquad
\lim_{|p|\to\infty} H^j(x,p) = +\infty \quad \textrm{unif. in $x$}.
\end{equation}

%
%
%
%
\section{Three definitions of viscosity solution}\label{sect:definition}

In this section, we recall the three definitions of viscosity solutions of problem~\eqref{HJ} introduced in \cite{acct}, \cite{imz} and \cite{sc}. Even though in \cite{imz} it is considered an evolutive equation, in order to compare the different notions of solution we restate the definition in terms of the  stationary equation \eqref{HJ}. Moreover
we will only consider the definitions at the vertex $v$, since in the other points of the network they coincide with the standard one of viscosity solution.

We first define the  admissible test functions  in the sense of \cite{acct} and  \cite{imz}
\begin{Definition}\label{accttest}
A function $\phi\in C(\overline{\G})$ is an $\{ACCT,IMZ\}$-admissible test function  at $v$ if, for any $j\in J$, $\phi^j$ belongs to $C^1([0,l_j])$.
We denote by $C^1_*(\G)$ the set of $\{ACCT,IMZ\}$-admissible test functions.
\end{Definition}
We now give the definition of admissible test function in the sense of \cite{sc}.
\begin{Definition}\label{cstest}
Let $\phi\in C(\G)$,  $j,k\in J$, $j\neq k$. A function $\phi:\G\to\R$ is said a $CS$-admissible $(j,k)$-test function at $v$, if
  $\phi$ is  differentiable at $ \pi_j^{-1}(v)$ and $ \pi_k^{-1}(v)$, respectively,  and
\begin{equation}\label{1:2}
D_j \phi(v  )+D_k \phi(v)=0.
\end{equation}

\end{Definition}
\begin{Remark}
In Definition \ref{accttest} admissible test functions can have different derivate  at $v$ along different edges, while in Definition \ref{cstest} an admissible $(j,k)$-test function is differentiable if restricted to couple $e_j$, $e_k$ taking into account the orientation.
On the other hand, In Definition \ref{accttest} admissible test function have derivatives along each incident arc, while in Definition \ref{cstest} an admissible test function needs to be differentiable only along two arcs.
\end{Remark}
%
%
\subsection{The ACCT solution}
In \cite{acct}, the Hamiltonian is the control-theoretic one
\begin{equation}\label{acctH}
    H(x,p)=\sup_{a\in A}\{-f(x,a)\cdot p-\ell(x,a)\}
\end{equation}
where $f$, $\ell$ are continuous functions while $A$ is a compact set.
It is assumed that $A=\cup_{j\in J} A^j$ where for $x\in e_j$,   $f(x,a)\in \R\eta_j$ if and only if $a\in A^j$.
In particular, we can write $H^j=\sup_{\a \in A^j}\{-f\cdot p-\ell\}$.
In \cite{acct}, the authors introduced a {\it relaxed gradient} for a function $u$ defined on $\G$. Consider $\zeta\in \R^2$ such that there exists a continuous function $z: [0,1]\to \G$ and a sequence $(t_n)_{n\in \N}$, $0 < t_n\le 1$ with $t_n\to 0$ such that
$\lim_{n\to \infty} \frac{z(t_n)-x}{t_n}=\zeta$, and such that  $\lim_{n\to \infty} \frac{u(z(t_n))-u(x)}{t_n}$ exists and does not depend on $z$ and $(t_n)_{n\in \N}$. In this case, they define the relaxed gradient by
\[Du(x,\zeta):=\lim_{n\to \infty} \frac{u(z(t_n))-u(x)}{t_n}.\]
For $\R_+:=[0,+\infty)$, we set
\[FL_j(v):= \overline{\rm co}\left ( (f(v,a),\ell(v,a)) : a\in A^j \right ) \cap (\R^+ \eta_j\times \R),\quad FL(v):=  \cup_{j\in J} FL_j(v). \]
\begin{Definition}\label{acctdef}
\begin{itemize}
\item A function $u\in USC(\G)$ is a subsolution of \eqref{HJ} at $v$
 if for any  $\phi\in C^1_*(\G)$ s.t. $u-\phi$ has a local maximum point at $v$, there holds
 \[
  u(v)+\sup_{ (\zeta,\xi) \in FL(v)}\{-D\phi(v,\zeta) -\xi\}\le 0.
 \]
\item  A function $u\in LSC(\G)$ is a  supersolution of \eqref{HJ}  at $v$ if for any   $\phi\in C^1_*(\G)$ s.t. $u-\phi$ has a local minimum point at $v$, there holds
\[
  u(v)+\sup_{ (\zeta,\xi) \in FL(v)}\{-D\phi(v,\zeta) -\xi\}\ge 0.
\]
\end{itemize}
\end{Definition}
\begin{Remark}
 Even if the coefficients in \eqref{acctH} are continuous functions, the Hamiltonian is in general discontinuous as a function
 of the state variable.  In fact,  \eqref{acctH} is the Hamiltonian of a control problem in $\R^2$ constrained to the network $\G$. Hence the
  set of the admissible controls, i.e. the controls corresponding to a tangential direction to the network, displays a discontinuity when passing from a point inside an edge to the vertex $v$.
 \end{Remark}
Let us now rewrite the previous definition in terms of the derivatives with respect to the parametrization:
\begin{Lemma}\label{acctdef2}
A function $u\in USC(\G)$ (respectively, $u\in LSC(\G)$)
 is a subsolution (resp., supersolution) of \eqref{HJ} at $v$ if for any $\phi\in C^1_*(\G)$ s.t. $u-\phi$ attains a local maximum (resp., minimum) at $v$, there holds
 \[
  u(v)+\max_{j\in J}\sup_{ (\zeta,\xi) \in FL_j(v)}\{-D_j \phi(v) \eta_j\cdot \zeta -\xi\}\le 0\qquad \textrm{(resp., $\geq 0$).}
 \]
\end{Lemma}

\subsection{The IMZ solution}
In \cite{imz}, the Hamiltonian $H$ is assumed to satisfy
 \begin{equation}\label{A1'}
 \text{$\exists p_0^j\in\R$ s.t.   $ H^j$ is non-increasing on $(-\infty,p_0^j)$,  non-decreasing on $(p_0^j,+\infty)$}
 \end{equation}
 for each $j\in J$. Define  the function $H^-_j$   by
 \begin{equation}\label{Hj-}
 H^-_j(p) := \inf_{q\le 0} H^j(v,p + q)= \left\{
                                        \begin{array}{ll}
                                          H^j(v,p), & \hbox{if $p<p_0^j$} \\
                                          H^j(v,p_0), & \hbox{if $p\ge p_0^j$.}
                                        \end{array}
                                      \right.
 \end{equation}
\begin{Definition}\label{imzdef}
\begin{itemize}
\item A function $u\in USC(\G)$ is a subsolution of \eqref{HJ} at $v$ if for any   $\phi\in C^1_*(\G)$ s.t. $u-\phi$ has a local maximum point at $v$, then
\[ u(v)+  \max_{j\in J} H^-_j( D_j\phi(v))\le 0.\]
\item A function $u\in LSC(\G)$ is a supersolution of \eqref{HJ}  at $v$ if for any $\phi\in C^1_*(\G)$ s.t. $u-\phi$ has a local minimum point at $v$, then
\[u(v)+\max_{j\in J} H^-_j( D_j\phi(v))\ge 0.\]
\end{itemize}
\end{Definition}
\begin{Remark}
In \cite{imz} the Hamiltonian is assumed to be independent of $x$ and strictly convex in $p$, but, for the purposes of the present paper, it suffices to require assumption \eqref{A1'}. It is important to observe that  no continuity condition on the Hamiltonian in $p$ at $v$ is assumed.
\end{Remark}

\subsection{The CS solution}
In \cite{sc}, the Hamiltonian is assumed  to satisfy
\begin{align}
 &H^j(v,p)=H^k(v,p)\quad\text{for any   $p\in \R$,  $j,k \in J$}\label{continuity}\\
 &H^j(v,p)=H^j(v,-p)\quad\text{for any  $p\in \R$,  $j\in J$.}\label{simmetry}
\end{align}
Assumptions~\eqref{continuity} and \eqref{simmetry} are the  continuity of $H$ in $p$ and its independence on the orientation of the incident arc, respectively.
\begin{Definition}\label{csdef}
\begin{itemize}
 \item A function $u\in USC(\G)$ is a subsolution of \eqref{HJ} at $v$ if
  for any $j,k\in J$ and any $(j,k)$-test function  $\phi$     for which
  $u-\phi$ attains a local maximum at $v$ relatively to $\bar e_j\cup \bar e_k$, then
  \begin{equation}\label{cssubsol}
u(v)+  H^k(v, D_k \phi(v) )=u(v)+  H^j(v, D_j \phi(v) )\le 0.
  \end{equation}
\item  A function $u\in LSC(\G)$ is a supersolution of \eqref{HJ}  at $v$ if  for any  $j\in J$, there exists $k\in J$, $k\neq j$, (said $v$-feasible for $j$ at $v$) such that for any  $(j,k)$-test function  $\phi$
   for which $u-\phi$ attains a local minimum at $v$ relatively to $\bar e_j\cup \bar e_k$, then
 \begin{equation}\label{cssupersol}
u(v)+  H^k(v, D_k \phi(v) )=u(v)+  H^j(v, D_j \phi(v) )\ge 0.
 \end{equation}
\end{itemize}
\end{Definition}
\begin{Remark}
Note that the    definitions of subsolution and supersolution in Definition \ref{csdef} are  not symmetric,  unlike  the ones in  Definitions \ref{acctdef} and \ref{imzdef}.
\end{Remark}
%
%
\section{Comparison between  ACCT and IMZ}\label{sect:acctimz}

To fix a common framework for the two settings, we assume that $H$ verifies \eqref{A1'} and it is given by \eqref{acctH}.

\begin{Theorem}\label{acctimz}
The definitions of   (ACCT)-solution and (IMZ)-solution  are equivalent.
\end{Theorem}
\begin{Proof}
We assume wlog $u(v)=0$  and that, for $p^-_j\leq \underline p_j<p_0^j<\bar p_j\leq p^+_j$, we have $H^j(v,p)< 0$ if, and only if, $p\in(\underline p_j,\bar p_j)$ and $H^j(v,p) >0$ if, and only if, $p\in (-\infty,p^-_j)\cup (p^+_j,+\infty)$.
Recall that the  admissible test functions coincide for the two definitions.\par
\noindent\texttt{1. (ACCT)-subsolution implies (IMZ)-subsolution.}
Let $u$ be an (ACCT)-subsolution and $\phi$ be an admissible upper test function for $u$ at $v$. To prove  $\max_j H^-_j ( D_j\phi(v))\leq 0$ we assume by contradiction  that $H^-_j( D_j\phi)>0$ for some $j\in J$. By the definition of $H^-_j$ this is equivalent to $D_j\phi(v)<p^-_j$. Hence, we have
\begin{align*}
\sup_{ (\zeta,\xi) \in FL(v)}\{-D\phi(v,\zeta) -\xi\}\ge
\sup_{(\zeta,\xi)\in FL_j(v)}\{-D_j\phi(v)\eta_j \cdot \zeta-\xi\}\\
 \geq \sup_{a\in A^j} \{-D_j\phi(v) \eta_j \cdot f(v,a)-\ell(v,a)\}\geq H^j(v,D_j\phi(v))>0
 \end{align*}
which contradicts the fact that $u$ is an (ACCT)-subsolution.\\
\noindent\texttt{2. (IMZ)-subsolution implies (ACCT)-subsolution.}
Let $u$ be a (IMZ)-sub\-so\-lu\-tion and $\phi$ be an admissible test function for $u$ at $v$. Assume  by contradiction
\begin{equation}\label{E100}
\sup_{(\zeta,\xi)\in FL(v)}\{-D\phi(v,\zeta)-\xi\}>0.
\end{equation}
Being a (IMZ)-subsolution, $u$ verifies $\max_j H^-_j(D\phi)\leq 0$ namely
\[D_j\phi(v)\geq p^-_j \qquad \forall j\in J.
\]
We recall that $FL_j(v)\subset \R_+ \eta_j\times \R$; therefore, the previous inequality implies
\[-D_j\phi(v)\cdot \zeta-\xi \leq -p^-_j\cdot \zeta-\xi \qquad \forall (\zeta,\xi)\in FL_j(v).\]
We deduce
\begin{equation}\label{E1}
\sup_{(\zeta,\xi)\in FL(v)}\{-D\phi (v, \zeta)-\xi\} \leq\max_j \sup_{(\zeta,\xi)\in FL_j(v)} \{-p^-_j \cdot \zeta-\xi\}.
\end{equation}
On the other hand, since  $H_j(v,p^-_j)=0$, by \eqref{acctH}  we obtain that
$-p^-_j\cdot f(v,a)-\ell(v,a)\leq 0$ for every $a \in A^j$. By linearity, we infer
$-p^-_j\cdot \zeta-\xi\leq 0$ for every $(\zeta,\xi)\in FL_j(v)$ and consequently
\[\max_j \sup_{(\zeta,\xi)\in FL_j(v)} \{-p^-_j \cdot \zeta-\xi\}\leq 0.\]
Replacing this inequality in \eqref{E1}, we get a contradiction to \eqref{E100}.\\
\noindent\texttt{3. (ACCT)-supersolution implies (IMZ)-supersolution.}
Let $u$ be a (ACCT)-supersolution. We want to prove that, for each admissible lower test function $\phi$  for $u$ at $v$, we have $\max_j H^-_j(D\phi(v))\geq 0$, i.e. that there exists $j\in J$ such that $D_j\phi\leq \underline p_j$.\\
We observe that $\pd_j^-u(v)=(-\infty,d_j]$ (possibly, $d_j=\pm\infty$) for each $j\in J$.
For $d_j=-\infty$, $\p_j^-u(v)$ is empty and, in particular, there is no lower test function for $u$ at $v$; thus, there is nothing to prove.
Assume $d_j\ne -\infty$ for every $j\in J$. We note that, for any admissible lower test function $\phi$, $D_j\phi(v)$ belongs to $(-\infty,d_j]$. If $d_j\leq \underline p_j$ for some $j\in J$, then there is nothing to prove.
By contradiction,   assume that $d_j>\underline p_j$ for every $j\in J$; hence, there exists an admissible lower test function $\phi$ such that $D_j\phi(v)\in(\underline p_j,\bar p_j)$ for each $j\in J$.
By Lemma \ref{acctdef2}, for some $j\in J$, we have
\begin{equation}\label{ac2}
\sup_{(\zeta,\xi)\in FL_j(v)}\{-D_j\phi(v)\cdot\zeta- \xi\}\geq 0.
\end{equation}
On the other hand, for each $j\in J$, there holds $H_j(v,D_j\phi(v))<0$ and, in particular,
\[
-D_j\phi(v)\cdot f(v,a)-\ell(v,a)<0\qquad \forall a\in A^j.
\]
By linearity, we infer
\[
-D_j\phi(v)\cdot \zeta-\xi< 0\qquad \forall (\zeta,\xi)\in FL_j(v)
\]
which contradicts~\eqref{ac2}.

\noindent\texttt{4. (IMZ)-supersolution implies (ACCT)-supersolution.} Let $u$ be a (IMZ)-supersolution and $\phi$ an admisible lower test function for $u$ at $v$. The definition of (IMZ)-supersolution reads $\max_j H^-_j(D_j\phi(v))\geq 0$, namely, for some $j\in J$, there holds $H^-_j(D_j\phi(v))\geq 0$.
This fact is equivalent to $D_j\phi(v)\leq \underline p_j$. Whence, we get $H_j(v,D_j\phi(v))=H^-_j(D_j\phi)\geq 0$ and we accomplish the proof observing that
$\sup_{(\zeta,\xi)\in FL(v)}\{-D\phi( v,\zeta)-\xi\}\geq \max_j H_j(v,D_j\phi(v))$.
\end{Proof}
\begin{Remark}
In the previous proof, we actually established that the definition of (ACCT)-supersolution (resp., subsolution) is equivalent to the one of (IMZ)-su\-per\-so\-lu\-tion (resp., subsolution).
\end{Remark}

%
%
\section{Comparison between IMZ and CS}\label{sect:imzcs}
We assume that $H$ satisfies \eqref{continuity}-\eqref{simmetry} and \eqref{A1'} with $p^j_0=0$ (because of \eqref{simmetry}).

\begin{Theorem}\label{csimz_super}
The definitions of (CS)-solution and (IMZ)-solution are equivalent.
\end{Theorem}

The proof of this Theorem is postponed at the end of the section. Let us first establish the following result.

\begin{Proposition}\label{imzcs}
The definitions of (IMZ)-subsolution and (CS)-subsolution  are equivalent, while (IMZ)-supersolution implies (CS)-supersolution.
\end{Proposition}
\begin{Proof}
For some positive values $\bar p$ and $p^*$, wlog we assume: $u(v)=0$, $H(v,p)>0$ if, and only if, $p\in (-\infty, -\bar p)\cup (\bar p,+\infty)$ and $H(v,p)<0$ if, and only if, $p\in (-p^*,p^*)$.

\noindent\texttt{1. (IMZ)-subsolution implies (CS)-subsolution.}
By \eqref{Hj-}, the function $H^-$ is defined as follows: $H^-(p)=H(v,p)$ if $p\leq 0$, $H^-(v,p)=H(0)<0$ if $p\geq 0$.
Let $u$ be a (IMZ)-subsolution. Following the same arguments of \cite[Lemma5.5]{imz}, one can easily prove that $u$ is Lipschitz continuous; in particular, $\pd^+_j u(v)\ne \emptyset$ for every $j\in J$. Set $\pd^+_j u(v)=:[p_j,+\infty)$ for each $j\in J$. We deduce that the function $\psi\in C^1_*(\G)$ with $\pd_j \psi=p_j$ is an (IMZ)-admissible upper test function. By the definition of (IMZ)-subsolution, we infer $H^- (p_j)\leq 0$ for each $j\in J$.
By the definition of $\bar p$ and of $H^-$, this relation can be rewritten as:
\begin{equation}\label{imzcs_st}
p_j\geq -\bar p \qquad \forall j\in J.
\end{equation}
Consider now a (CS)-admissible $(j,k)$-upper test function $\phi$ for $u$ at $v$. We want to prove that $H(v,D_j \phi)=H(v,D_k \phi)\leq 0$, namely $D_j \phi(v), D_k \phi(v)\in [-\bar p,\bar p]$.
To this end, we assume by contradiction that $D_j \phi(v)>\bar p$. Hence, by relation~\eqref{1:2}, we have
\begin{equation*}
-\bar p>-D_j\phi (v)=D_k \phi(v)\geq p_k
\end{equation*}
that contradicts relation~\eqref{imzcs_st}. Therefore, we have $D_j \phi(v)\leq \bar p$ and, similarly, $D_k \phi(v)\leq \bar p$. Owing to the relation $D_k \phi(v)=-D_j \phi(v)$, it  holds $D_j \phi(v), D_k \phi(v)\in [-\bar p,\bar p]$.
Taking into account the arbitrariness of the function~$\phi$ and of $(j,k)$, we accomplish the proof of point $(1)$.\\
\noindent\texttt{2. (CS)-subsolution implies (IMZ)-subsolution.}
We will use the following Lemma (see \cite[Lemma 5.4]{cm})
\begin{Lemma}\label{Lemma54}
Let $y_m\in [0,l_j]$ ($m\in \N$) with $\lim_m y_m=0$. Then there holds
 \begin{equation}\label{prp:noconvex}
 \lim_{m\to +\infty} H\left(v, \frac{u^j(y_m)}{y_m}\right)\le 0.
\end{equation}
\end{Lemma}
Assume by contradiction that there exists $\phi\in C^1_*(\G)$ and $j\in J$ such that $H^-_j(D_j\phi(v))\ge \d>0$. Set $p=D_j\phi(v)$ and fix $\eta$ sufficiently small  in such a way that $H^-_j(p+\eta)\ge \d/2$.
Since $p\in \pd_j^+u(v)$,  if $x_m\in e_j$, $x_m\to v$ and $y_m=\pi_j^{-1}(x_m)$ we have that
\[ p_m:=\frac{u^j(y_m)-u^j(0)}{y_m}\le p+\eta  \]
for $m$ sufficiently large. Moreover $u(x_m)\to u(v)=0$ and therefore for $m$ sufficiently large
$u(x_m)+ H^-_j(p_m)\ge \frac{\d}{4}$. On the other side by Lemma \ref{Lemma54} and $u(v)=0$, it follows  that
$u(x_m)+H(v,p_m)<\d/4$ for $m$ sufficiently large and therefore a contradiction.\\
\noindent\texttt{3. (IMZ)-supersolution implies (CS)-supersolution.}
Let us observe that it is possible that, for some $j\in J$, there holds $\pd^-_j u(v)= \emptyset$. In this case there is no (IMZ)-admissible lower test function. In order to prove that $u$ is a (CS)-supersolution, it suffices to choose the edge~$j$ as the feasible one; with this choice, there is no (CS)-admissible test function neither, i.e. there is nothing to prove.

Let us now assume $\pd^-_j u(v)\ne \emptyset$ for every $j\in J$. Set $\pd^-_j u(v)=:(-\infty,p_j]$ (possibly $p_j=+\infty$). We deduce that the function $\psi\in C^1_*(\G)$ with $D_j \psi(v)=p_j$ ($D_j \psi=2p^*$ if $p_j=+\infty$) is an (IMZ)-admissible lower test function. By the definition of (IMZ)-supersolution, there exists $j\in J$ such that $H^- (D_j\psi(v))\geq 0$. By the definition of $H^-$ and   of $\psi$ when $p_j=+\infty$, we infer that there exists $j\in J$ such that
\begin{equation}\label{imzcs_sp}
p_j\leq -p^*.
\end{equation}
In order to prove that $u$ is a (CS)-supersolution, let us choose $e_j$ as the feasible edge; in other words, we claim that, for every $k\in J\setminus \{j\}$, for every (CS)-admissible $(j,k)$-lower test function~$\phi$, there holds $H(v,D_j\phi(v))=H(v,D_k\phi(v))\geq 0$.
To this end, we observe that the definition of lower test function and relation~\eqref{imzcs_sp} entail: $D_j\phi(v)\leq p_j\leq -\bar p$. By \eqref{A1'}, we get $H(D_j\phi(v))\geq 0$.
\end{Proof}
\begin{Remark}
Thanks to Proposition \ref{acctimz}, we also have that the definitions of (ACCT)-subsolution  and  (CS)-subsolution  are equivalent, while (ACCT)-supersolution implies (CS)-supersolution.
\end{Remark}


\begin{Proofc}{Proof of Theorem \ref{csimz_super}} By Proposition \ref{imzcs}, we have only to prove that a (CS)-solution is a (IMZ)-supersolution.
Let $u$ be a (CS)-solution; in particular, let us recall that $u$ is Lipschitz continuous. We observe that, for each $j\in J$, there holds: $\pd_j^- u(v)=(-\infty,p_j]$ for some $p_j\in\R$.
A function~$\psi\in C^1_*(\G)$ is a (IMZ)-admissible lower test function if, and only if, $D_j\psi(v)\leq p_j$ for each $j\in J$. For such a $\psi$, our aim is to prove:
\begin{equation}\label{aa1}
u(v)+ H_j^-(D_j \psi)\geq0\qquad \textrm{for some }j\in J.
\end{equation}
Let us split our arguments according to the existence or not existence of (CS)-admissible lower test functions.
Fix an edge, say $e_1$; wlog assume that $e_2$ is the feasible edge in the definition of (CS)-supersolution.
We observe that there exists an admissible $(1,2)$-lower test function if, and only if, there holds $p_1\geq -p_2$.

\noindent\texttt{Case (i): Non existence of a (CS)-lower test function.}
We have
\begin{equation}\label{ab1}
p_1<-p_2.
\end{equation}
By \eqref{ab1} either $p_1$ or $p_2$ is negative and wlog we assume $p_1<0$. We observe that, for the network $\bar e_1 \cup \bar e_2$, the definition of (CS)-solution coincides with the classical definition of solution in viscosity sense for a segment. We can consider such a segment as the one given by the orientation of $e_1$   inverting the  orientation of $e_2$.
 By \cite[Theorem1]{JS} it follows that that either $u$ is differentiable in $v$ (this cannot be our case because $p_1\ne p_2$) or there holds
\begin{equation*}
H(v,p_1)=H(v,-p_2)=-u(v).
\end{equation*}
Since $p_1<0$, we deduce
\begin{equation*}
u(v)+H^- (p_1)=0,
\end{equation*}
and also, by the monotonicity of $H^-$,
\begin{equation*}
u(v)+H^-(p)\geq0, \qquad \forall p\leq p_1
\end{equation*}
which entails inequality~\eqref{aa1}.

\noindent\texttt{Case (ii): Existence of some (CS)-lower test function.}
We assume that there exists some (CS)-admissible $(1,2)$-lower test function~$\phi$ for the function~$u$ at $v$. In this case, we have $p_1\geq D_1\phi(v)=-D_2\phi(v)\geq -p_2$. Moreover, observe that any $(1,2)$-lower test function~$\bar \phi$ with $D_1\bar  \phi(v)\in[-p_2,p_1]$ is admissible. Therefore, the definition of (CS)-supersolution yields
\begin{eqnarray}
u(v)+H(v,p)\geq0, &\qquad& \forall p \in[-p_2,p_1], \label{ab2}\\
u(v)+H(v,p)\geq0, &\qquad& \forall p \in[-p_1,p_2].\label{ab3}
\end{eqnarray}
By \eqref{A1'}, the function $h(p):=u(v)+H(v,p)$ is negative if, and only if, $p\in(- p^*,p^*)$ for some positive value $p^*$. Relation~\eqref{ab2} ensures that the interval $(-p_2,p_1)\cap(-p^*,p^*)$ must be empty. Whence, we have either $p_1\leq -p^*$ or $p_2\leq -p^*$.
For $p_1\leq -p^*$, the definition of $H_1^-$ and relation~\eqref{ab2} ensure
\begin{equation*}
u(v)+H_1^-(p)=u(v)+H(v,p)\geq0, \qquad \forall p\leq p_1.
\end{equation*}
%
For $p_2\leq -p^*$, the definition of $H^-_2$ and relation~\eqref{ab3} ensure
\begin{equation*}
u(v)+H^-_2(p)=u(v)+H(v,p)\geq0, \qquad \forall p\leq p_2.
\end{equation*}
Both the last two relations amount to inequality~\eqref{aa1}.
\end{Proofc}

\begin{Remark}
In the previous proof the hypothesis that $u$ is a (CS)-solution is needed only when there is no admissible (CS)-lower test function. In this case, it seems that the information  available from the (CS)-definition of supersolution are not sufficient  to obtain an equivalent property for (IMZ)-definition.
\end{Remark}
\textbf{\large Acknowledgments: }
The authors wish to thank Yves Achdou and Nicoletta Tchou for many useful discussions on the subject.


\end{document}